\documentclass[11pt]{article}
\usepackage{amsfonts}
\usepackage{eufrak}
\usepackage{mathrsfs}
\usepackage{amssymb}
\usepackage{mathtools}
\usepackage{commath}
\usepackage{dsfont}
\begin{document}
\title{\bf On a binary Diophantine inequality\\ involving primes of a special type
\thanks{Project supported by the National Natural Science Foundation of China (grant no.11771333) and by the Fundamental Research Funds for the Central Universities (grant no. JUSRP122031). }}
\author{{Yuhui  Liu}\\ {School of Science,  Jiangnan University}\\{Wuxi 214122, Jiangsu, China}
\\{Email: tjliuyuhui@outlook.com}}
\date{}
\maketitle
{\bf Abstract}
Let $1<c<\frac{1787}{1502}$ and $N$ be a sufficiently large real number. In this paper, it is proved that for any arbitrarily large number $E>0$ and for almost all real $R \in (N,2N]$, the Diophantine inequality
\begin{align*}
\left|p_{1}^{c}+p_{2}^{c}-R\right|<\left(\log N\right)^{-E}
\end{align*}
is solvable in prime variables $p_1,p_2$ such that, each of the numbers $p_{1}+2,\,p_{2}+2$ has at most $\big[\frac{79606}{35740-30040c}\big]$ prime factors, counted with multiplicity. Moreover, we prove that the Diophantine inequality
\begin{align*}
\left|p_{1}^{c}+p_{2}^{c}+p_{3}^{c}+p_{4}^{c}-N\right|<\left(\log N\right)^{-E}
\end{align*}
is solvable in prime variables $p_1,p_2,p_3,p_4$ such that, each of the numbers $p_{i}+2\,(i=1,2,3,4)$ has at most $\big[\frac{93801402}{35740000-30040000c}\big]$ prime factors, counted with multiplicity.



\section{Introduction}
\setcounter{equation}{0}
\hspace{5.8mm} Assume that $c>1,c\not\in\mathbb{N}$ and let $\Xi$ be a small positive number depending on $c$. For sufficiently large real number $N$, let $H(c)$ denote the smallest positive integer $r$ such that the inequality
\begin{align}
\left|p_{1}^{c}+p_{2}^{c}+\cdots+p_{r}^{c}-N\right|<\Xi
\end{align}
is solvable in prime numbers $p_1,p_2,\cdots,p_r$.

In 1952, Piatetski-Shapiro [16] first proved that
\begin{align*}
\lim\sup\limits_{c\rightarrow\infty} \frac{H(c)}{c \log c} \le 4.
\end{align*}
He also proved that if $1<c<\frac32$, then $H(c)\leq 5.$ From this result and the Goldbach-Vinogradov Theorem, it is reasonable to conjecture that if $c$ is close to $1$, one should expect $H(c) \le 3$.

In 1999, Laporta [10] considered an analogous result for the corresponding binary problem. He proved that for any sufficiently large real number $N$ and for almost all real $R \in (N,2N]$ (in the sense of Lebesgue measure), the inequality
\begin{align}
\left|p_{1}^{c}+p_{2}^{c}-R\right|<\varepsilon
\end{align}
is solvable in prime variables $p_1,p_2$, where $1<c<\frac{15}{14}$ and $\varepsilon$ relies on $R$. Later, the upper bound of $c$ was enlarged to
\begin{align*}
\frac{43}{36}=1.194...\,[21],\,\,\,\,\,\,\frac65\,[5],\,\,\,\,\,\,\frac{59}{44}=1.341...\,[14],\,\,\,\,\,\,\frac{39}{29}=1.345...\,[1].
\end{align*}

On the other hand, let $\mathcal{P}_{r}$ be an almost-prime with at most $r$ prime factors, counted according to multiplicity. One well-known result is due to Chen [3,4], who proved that there exist infinitely many primes $p$ such that $p+2 \in \mathcal{P}_{2}$. Motivated by Chen [3,4], it is reasonable to conjecture that for some constants $c$, the inequality (1.1)  with a suitable $\Xi:=\:\Xi(N)$ satisfying $\Xi(N)\rightarrow0$ as $N\rightarrow\infty$, is solvable in primes $p_i$ such that $p_{i}+2 \in \mathcal{P}_{r}$ for a fixed $r\ge2$.

In 2017, Tolev [18] considered the case $s=3$ and proved that for $1<c<\frac{15}{14}$, the Diophantine inequality (1.1) with $s=3$ and $\Xi=:(\log N)^{-E}$, where $E>0$ is an arbitrarily large constant, has a solution in primes $p_{1},\:p_{2},\:p_{3}$ such that each of the numbers $p_{i}+2$ for $i\:=\:1,2,3$ has at most $[\frac{369}{180-168c}]$ prime factors, counted with multiplicity. Very recently, the bound for $c$ was improved to $\frac{281}{250}$ by Zhu [22], and then to $\frac{973}{856}$ in Li, Xue and Zhang [11], while the number of prime factors that each of the numbers $p_{i}+2$ for $i=1,2,3$ has is $\left[{\frac{1475}{562-500c}}\right]$ [22] and $\left[{\frac{12626}{4865-4280c}}\right]$ [11], counted with multiplicity.

For the case $s=2$ for this  type, motivated by Li, Cai [14], Baker [1] and Chen [3,4], it is reasonable to conjecture that if $E>0$ be an arbitrarily large constant and $c$ is close to $1$, then for almost all real $R \in (N,2N]$, the inequality (1.2) with $\varepsilon=:(\log N)^{-E}$ is solvable in primes $p_1,p_2$ such that $p_1+2$ and $p_2+2$  are almost-primes of a certain fixed order. In this paper, we investigate the above problem by giving the following Theorem 1.
\\
\\
{\bf Theorem 1} {\it \,\,If $1<c<\frac{1787}{1502}$ and $N$ is a sufficiently large real number, then for any arbitrarily large numbers $E>0$ and $A>0$, for all real $R \in (N,2N] \setminus \mathcal{A}$ with $|\mathcal{A}|=O\bigg(N\left(\frac{\log N}{c}\right)^{\frac92-A}\bigg),$ the Diophantine inequality
\begin{align}
\left|p_{1}^{c}+p_{2}^{c}-R\right|<\left(\log N\right)^{-E}
\end{align}
is solvable in prime variables $p_1,p_2$ such that each of the numbers $p_{i}+2\,\, (i=1,2)$ has at most $\big[\frac{79606}{35740-30040c}\big]$  prime factors, counted with multiplicity, where $|\mathcal{A}|$ denotes the cardinality of the set $\mathcal{A}$.}

In 2003, Zhai and Cao [20] showed that the inequality (1.1) is solvable for $r=4$  in the range $1<c<\frac{81}{68}$. Later, the range was enlarged to $\frac{97}{81}$,\,$\frac65$,\,$\frac{59}{44}$,\,$\frac{1198}{889}$, by Liu and Shi [15], Li and Zhang [12], Li and Cai [14], Li and Zhang [13] successively. The best result till now is due to Baker [1], and he enlarged the range to $1<c<\frac{39}{29}$. Motivated by Baker [1],  Li and Cai [14]  and Theorem 1, we have the following:
\\
\\
{\bf Theorem 2} {\it \,\,If $1<c<\frac{1787}{1502}$ and $N$ is a sufficiently large real number, then for any arbitrarily large number $E>0,$ the Diophantine inequality
\begin{align}
\left|p_{1}^{c}+p_{2}^{c}+p_{3}^{c}+p_{4}^{c}-N\right|<\left(\log N\right)^{-E}
\end{align}
has a solution in prime numbers $p_1,p_2,p_3,p_4$ such that, for each $i\in \{1,2,3,4\}$,  each of the numbers $p_{i}+2$ has at most $\big[\frac{93801402}{35740000-30040000c}\big]$ prime factors, counted with multiplicity.}

\section{Notation and some preliminary lemmas}
\setcounter{equation}{0}
\hspace{5.8mm}For the proof of the Theorems, in this section we introduce the necessary notation and lemmas.

Throughout this paper,  by $A$ and $E$ we denote sufficiently large constants. In addition, let $\varepsilon < 10^{-100}$ be an arbitrarily small positive constant not necessarily the same in different formulae.  The letter $p$, with or without subscripts, is reserved for a prime number. We use $e(\alpha)$ to denote $e^{2\pi i\alpha}$. By $A \sim B$ we mean that $B < A \leqslant 2B$ and $f(x) \asymp g(x)$ means that $f(x) \ll g(x) \ll f(x)$. We denote by $(m,n)$ the greatest common divisor of $m$ and $n$. As usual, we use $\varphi(n)$, $\mu(n)$, $\Lambda(n)$ to denote Euler's function, M\"{o}bius' function and von Mangoldt's function. Let
\begin{eqnarray*}
&&1<c<\frac{1787}{1502},\,\,X=N^{\frac1c},\,\,\delta=\frac{1787}{1502}-c,\,\,\xi=\frac{2c}{3}-\frac{1}{3},\,\,\eta_1=\frac{20}{53}\delta,\,\,
\eta_2=\frac{20000}{62451}\delta,\,\,\mathcal{L}=\log X,\\
&&\mbox{}\nonumber\\
&&\Delta=(\log N)^{-E},\,\,D=X^{\delta},\,\,\tau=X^{\xi-c},\,\,P(z)=\prod_{2<p<z_j}p,\,\,z_j=X^{\eta_j}(j=1,2),\,\,K=\mathcal{L}^{E+3}.
\end{eqnarray*}
For fixed constant $\mu\in(0,1)$, $\lambda^{\pm}(d)$ be the Rosser functions of level $D$ and $\lambda(d)$ be real numbers satisfying
\begin{align*}
|\lambda(d)|\leq1,\quad\lambda(d)=0\quad\mathrm{if}\quad2|d\quad\mathrm{or}\quad\mu(d)=0.
\end{align*}
For $i=1,2,3,4$, we denote
\begin{align}
\Lambda_i=\sum_{d|(p_i+2,P(z))}\mu(d)=\begin{cases}1,&\mathrm{if}\left(p_i+2,P(z)\right)=1,\\0,&\mathrm{otherwise},\end{cases}
\end{align}
and
\begin{align}
\Lambda_{i}^{\pm}=\sum_{d|(p_{i}+2,P(z))}\lambda^{\pm}(d).
\end{align}
Moreover, we define
\begin{eqnarray*}
&&I(x)=\int_{\mu X}^{X}e(t^{c}x)\mathrm{d}t,\\
&&\mbox{}\nonumber\\
&&L(x)=\sum_{d\leq D}\lambda(d)\sum_{\mu X<p\leq X\atop{d|p+2}}(\log p)e(p^cx),\\
&&L^\pm(x)=\sum_{d|P(z)}\lambda^\pm(d)\sum_{\mu X<p\leq X\atop{d|{p+2}}}(\log p)e(p^cx).
\end{eqnarray*}
Next, we state the lemmas required in this paper.\\

{\bf Lemma 2.1}. {\it There exists a function $\varphi(y)$ which is $r = [\mathcal{L}^2]$ times continuously differentiable and satisfies
\begin{align*}
\varphi(y)& =1\quad for\,\,\left|y\right|\leq a-\Delta,  \\
0<\varphi(y)& <1\quad for\,\,\,a-\Delta<|y|<a+\Delta,  \\
\varphi(y)& =0\quad for\,\,|y|\geq a+\Delta,
\end{align*}
and its Fourier transform
\begin{align*}
\Theta(x)=\int_{-\infty}^{\infty}e(-xy)\varphi(y)\mathrm{d}y
\end{align*}
satisfies the inequality}
\begin{align*}
|\Theta(x)|\leq\min\left(2a,\frac{1}{\pi|x|},\frac{1}{\pi|x|}\left(\frac{r}{2\pi|x|\Delta}\right)^{r}\right).
\end{align*}

{\bf Proof}. See Piatetski-Shapiro [16].\\

{\bf Lemma 2.2}. {\it Suppose that $D>4$ is a real number. There exist Rosser functions $\lambda^{\pm}(d)$ of level $D$ with the following properties:}\\
\mbox{(i)} {\it For any positive integer $d$ and $n$, we have
\begin{align*}
|\lambda^{\pm}(d)|\leq1,\quad\lambda^{\pm}(d)=0\,\, if\,\, d\ge D\,\, or\,\, \mu(d)=0,
\end{align*}
and}
\begin{align*}
\sum_{d|n}\lambda^{-}(d)\leq\sum_{d|n}\mu(d)\leq\sum_{d|n}\lambda^{+}(d).
\end{align*}
\mbox{(ii)} {\it For real number $z$, define
\begin{align}
\mathfrak{P}=\prod_{2<p<z}\left(1-\frac{1}{p-1}\right),\quad\mathcal{M}^{\pm}=\sum_{d|P(z)}\:\frac{\lambda^{\pm}(d)}{\varphi(d)},\quad s_{0}=\frac{\log D}{\log z}.
\end{align}
Then we have
\begin{align}
&\mathfrak{P} \asymp {\mathcal{L}}^{-1},\,\,\,\,\,\,\,\,\mathcal{M}^{\pm}\ll\sum_{d\leq D}\frac1{\varphi(d)}\ll \mathcal{L},\\
&\mathfrak{P}\le\mathcal{M}^+\le\mathfrak{P}\Big(F(s_0)+O\big((\log D)^{-1/3}\big)\Big),\notag\\
&\mathfrak{P}\ge\mathcal{M}^-\ge\mathfrak{P}\Big(f(s_0)+O\big((\log D)^{-1/3}\big)\Big),\notag
\end{align}
where $F\left(s\right)$ and $f\left(s\right)$ are defined by
\begin{align*}
F(s)=\begin{cases}2e^\gamma s^{-1},&if~1\leq s\leq 3,\\
2e^\gamma s^{-1}\left(1+\int_2^{s-1}t^{-1}\log(t-1)\mathrm{d}t\right),&if~3\leq s\leq 5,
\end{cases}
\end{align*}
and
\begin{align*}
f(s)=2e^\gamma s^{-1}\log(s-1),\quad if\:2\leq s\leq4,
\end{align*}
where $\gamma$ denotes the Euler's constant}.\\

{\bf Proof}. We refer the reader to Greaves [6, \S 4].\\

{\bf Lemma 2.3}. {\it Suppose that $\Lambda_{i}$ and $\Lambda_{i}^{\pm} \,(i=1,2,3,4)$ are real numbers defined by} (2.1) {\it and} (2.2). {\it Then we have}
\begin{align*}
\mbox{(i)}& \,\,\Lambda_{1}\Lambda_{2}\geq\Lambda_{1}^{-}\Lambda_{2}^{+}+\Lambda_{1}^{+}\Lambda_{2}^{-}-\Lambda_{1}^{+}\Lambda_{2}^{+};\\
\mbox{(ii)}& \,\,\Lambda_{1}\Lambda_{2}\Lambda_{3}\Lambda_{4}\geq\Lambda_{1}^{-}\Lambda_{2}^{+}\Lambda_{3}^{+}\Lambda_{4}^{+}+\Lambda_{1}^{+}\Lambda_{2}^{-}\Lambda_{3}^{+}\Lambda_{4}^{+}
+\Lambda_{1}^{+}\Lambda_{2}^{+}\Lambda_{3}^{-}\Lambda_{4}^{+}+\Lambda_{1}^{+}\Lambda_{2}^{+}\Lambda_{3}^{+}\Lambda_{4}^{-}-3\Lambda_{1}^{+}\Lambda_{2}^{+}\Lambda_{3}^{+}\Lambda_{4}^{+}.
\end{align*}

{\bf Proof}. The proof is elementary and similar to the proof of J. Br\"{u}dern and E. Fouvry [2, Lemma 13].\\

{\bf Lemma 2.4}. {\it Suppose that $f(x):[a,b]\rightarrow\mathbb{R}$ has continuous derivatives of arbitrary order on $[a,b]$, where $1\leqslant a<b\leqslant 2a$. Suppose further that
\begin{align*}
\left|f^{(j)}(x)\right|\asymp\lambda_{1}a^{1-j},\quad j\geqslant1,\quad x\in[a,b].
\end{align*}
Then for any exponential pair $(\kappa,\lambda),$ we have}
\begin{align*}
\sum_{a<n\leqslant b}e(f(n))\ll\lambda_{1}^{\kappa}a^{\lambda}+\lambda_{1}^{-1}.
\end{align*}

{\bf Proof}. See (3.3.4) of Graham and Kolesnik [5].\\

{\bf Lemma 2.5}. {\it $(\frac{187}{659}+\varepsilon,\,\frac{374}{659}+\varepsilon)$ is an exponential pair}.

{\bf Proof}. From Huxley [8, Table 17.3], $(\frac{89}{570}+\varepsilon,\,\frac{374}{570}+\varepsilon)$ is an exponential pair. Then $BA(\frac{89}{570}+\varepsilon,\,\frac{374}{570}+\varepsilon)= (\frac{187}{659}+\varepsilon,\,\frac{374}{659}+\varepsilon)$ is an exponential pair.\\

{\bf Lemma 2.6}. {\it Suppose that $\chi$ is a primitive Dirichlet character modulo $d$. For $Q\geq1,T\geq2$, let
\begin{align*}
\mathcal{L}(T,Q,X)=\sum_{d\leq Q}\sum_{\chi(d)^{*}}\sum_{|\gamma(\chi)|\leq T}X^{\beta(\chi)},
\end{align*}
where the summation in the inner most sum is taken over the non-trivial zeros $\rho(\chi)=\beta(\chi)+i\gamma(\chi)$ of Dirichlet's L-function $L(s,\chi)$ such that $|\gamma(\chi)|\leq T$ . Then we have}
\begin{align*}
&\mbox{(i)}\,\, \sum_{d\leq Q}\sum_{\chi(d)^{*}}N(T,\sigma,\chi)\ll(Q^{2}T)^{\frac{12(1-\sigma)}{5}}(\log(QT))^{44},\\
&\mbox{(ii)}\,\mathcal{L}(T,Q,X)\ll X{\mathcal{L}}^{45}+X^{\frac{1}{2}}Q^{\frac{12}{5}}T^{\frac{6}{5}}{\mathcal{L}}^{45}.
\end{align*}

{\bf Proof}. For (i), see Hulexy [7]. For (ii), it follows from the same argument as Lemma 4.3 in Zhu [22] excepet that we use Lemma 2.6 (i) instead of Lemma 4.3 in Zhu [22].\\

{\bf Lemma 2.7}. {\it If $2\xi + 16\delta <5$ and $7\xi + 14\delta <5$, then for $|x|\leq\tau,$ we have}
\begin{align*}
&\mbox{(i)}\,\, L(x) =\sum_{d\leq D}\frac{\lambda(d)}{\varphi(d)}I(x)+ O\big(X{\mathcal{L}}^{-A}\big),\\
&\mbox{(ii)}\,L^{\pm}(x)=\mathcal{M}^{\pm}I(x)+ O\big(X{\mathcal{L}}^{-A}\big),
\end{align*}
{\it where $\mathcal{M}^{\pm}$ is defined by} (2.3).

{\bf Proof}. The proof is exactly the same as Lemma 4.5 in  Zhu [22] by using Lemma 2.6 (ii) instead of Lemma 4.4 in   Zhu [22].\\

{\bf Lemma 2.8}. {\it We have}
\begin{align*}
\mbox{(i)}&\int_{|x|\leq\tau}\left|L(x)\right|^2\mathrm{d}x\ll X^{2-c}{\mathcal{L}}^{6};\,\,\,\,\,\,\,\,\,\,\,\mbox{(ii)}\int_{|x|\leq\tau}\left|L(x)\right|^3\mathrm{d}x\ll X^{3-c}{\mathcal{L}}^{8};\\
\mbox{(iii)}&\int_{|x|\leq\tau}\left|I(x)\right|^2\mathrm{d}x\ll X^{2-c}{\mathcal{L}}^4;\,\,\,\,\,\,\,\,\,\,\,\mbox{(iv)}\int_{|x|\leq\tau}\left|I(x)\right|^3\mathrm{d}x\ll X^{3-c}{\mathcal{L}}^4;\\
\mbox{(v)}&\int_{|x|\leq K}\left|L(x)\right|^3\mathrm{d}x\ll X^2K{\mathcal{L}}^8.
\end{align*}

{\bf Proof}. By using the trivial estimate $|L(x)|\ll X{\mathcal{L}}^2$ and $|I(x)|\ll X$ , the conclusions follow from Lemma 11 of Tolev [18].\\

{\bf Lemma 2.9}. {\it We have}
\begin{align*}
\mbox{(i)}&\int_{-\infty}^{+\infty} I^2(x)e(-xR)\Theta(x)\mathrm{d}x\gg \Delta R^{\frac{2}{c}-1};\\
\mbox{(ii)}&\int_{-\infty}^{+\infty} I^4(x)e(-xN)\Theta(x)\mathrm{d}x\gg \Delta X^{4-c}.
\end{align*}

{\bf Proof}. See Kumchev and Laporta [9] and Zhai and Cao [20] respectively.\\

{\bf Lemma 2.10}. {\it For $1<c<\frac{1787}{1502}$, we have}
\begin{align*}
\sup_{\tau\leq|x|\leq K}\left|L(x)\right|\ll X^{\frac{4}{3}-\frac{1}{3}c-\varepsilon}.
\end{align*}

{\bf Proof}. We follow the notations and argument as in Tolev [18]. Let
\begin{align*}
L_1(x)= \sum_{d\leq D}\lambda(d)\sum_{\mu X<n\leq X\atop{d|n+2}}\Lambda(n)e(n^cx).
\end{align*}
It is obvious that
\begin{align*}
L(x)= L_1(x)+O\big(X^{\frac{1}{2}+\varepsilon}\big).
\end{align*}
So in the following, we consider the sum $L_1(x)$ instead of $L(x)$. We rewrite $L_{1}(x)$ in the form
\begin{align*}
L_{1}(x)=\sum_{\mu X<n\leq X}\Lambda(n)f(n),\,\,\,\,\,f(n)=\sum_{d\leq D\atop{d|n+2}}\lambda(d)e(n^{c}x).
\end{align*}
By the arguments in Vaughan [19,\,p.112], we further divide $L_1(x)$ into three parts:
\begin{align*}
L_{1}(x)=S_{1}-S_{2}-S_{3},
\end{align*}
where\\
\begin{align*}
S_1 \,&=\sum_{k\leq X^{1/3}}\mu(k)\sum_{\frac{\mu X}{k}<\ell\leq\frac{X}{k}}(\log\ell)f(k\ell),  \\
S_2 \,&=\sum_{k\leq X^{2/3}}c(k)\sum_{\frac{\mu X}{k}<\ell\leq\frac{X}{k}}f(k\ell) \\
& = \sum_{k\leq X^{1/3}}c(k)\sum_{\frac{\mu X}{k}<\ell\leq\frac{X}{k}}f(k\ell) + \sum_{X^{1/3}<k\leq X^{2/3}}c(k)\sum_{\frac{\mu X}{k}<\ell\leq\frac{X}{k}}f(k\ell)\\
& =: S_{2}^{'} + S_{2}^{''},\\
S_3 \,&=\sum_{X^{1/3}<k\leq X^{2/3}}a(k)\sum_{\frac{\mu X}{k}<\ell\leq\frac{X}{k}}\Lambda(\ell)f(k\ell),
\end{align*}
and $|a(k)|\leq \tau(k),|c(k)|\leq\log k$. Hence
\begin{align}
L_{1}(x)\ll\left|S_{1}\right|+|S_{2}^{'}|+|S_{2}^{''}|+\left|S_{3}\right|.
\end{align}

First, we consider the sum $S_{2}^{'}$. From the argument (4.12)-(4.13) in Tolev [18], we obtain
\begin{align}
S_{2}^{'}=\sum_{d\leq D}\lambda(d)\sum_{k\leq X^{1/3}\atop{(k,d)=1}}c(k)\sum_{\frac{\mu X}{kd}-\frac{\ell_{0}}{d}<m\leq\frac{X}{kd}-\frac{\ell_{0}}{d}}e\big(h(m)\big),\,\,\,\,\,\,\,h(m)=xk^{c}(\ell_{0}+md)^{c}.
\end{align}
Since
\begin{align*}
\left|h^{(j)}(m)\right|\asymp\left(|x|dkX^{c-1}\right)\cdot\left(Xk^{-1}d^{-1}\right)^{1-j},\quad j\geq1.
\end{align*}
Using Lemma 2.4 with $(\kappa,\lambda)=ABA^3B(0,1)=(\frac{11}{82},\frac{57}{82})$, we know that the inner sum over $m$ in (2.6) is
\begin{align}
&\ll\left(|x|dkX^{c-1}\right)^{\frac{11}{82}}\left(Xk^{-1}d^{-1}\right)^{\frac{57}{82}}+\left(|x|dkX^{c-1}\right)^{-1}\notag\\
&\ll|x|^{\frac{11}{82}}X^{\frac{11}{82}c+\frac{23}{41}}k^{-\frac{23}{41}}d^{-\frac{57}{82}}+|x|^{-1}k^{-1}d^{-1}X^{1-c}.
\end{align}
Inserting (2.7) into (2.6), we find that
\begin{align}
S_{2}^{'}\ll X^{\varepsilon}\left(D^{\frac{25}{82}}|x|^{\frac{11}{82}}X^{\frac{11c}{82}+\frac{29}{41}}+|x|^{-1}X^{1-c}\right).
\end{align}
We apply Abel's summation formula to get rid of the factor $\log\ell$ in the definition of $S_1$ and then proceed as the process of $S_{2}^{'}$ to treat the sum $S_1$. In this way, we find that
\begin{align}
S_{1}\ll X^{\varepsilon}\left(D^{\frac{25}{82}}|x|^{\frac{11}{82}}X^{\frac{11c}{82}+\frac{29}{41}}+|x|^{-1}X^{1-c}\right).
\end{align}

Consider now the sum $S_3$. By a splitting argument, we divide $S_3$ into $O(\mathcal{L})$ sums of the type
\begin{align*}
W(K):=\sum_{K<k\leq K_1}a(k)\sum_{\frac{\mu X}{k}<\ell\leq\frac{X}{k}}\Lambda(\ell)\sum_{d\leq D\atop{d|k\ell+2}}\lambda(d)e\big((k\ell)^cx\big),
\end{align*}
where
\begin{align*}
K<K_{1}\leq 2K,\quad X^{\frac{1}{3}}\leq K\:<K_{1}\leq X^{\frac{2}{3}}.
\end{align*}
We first consider the case $K\geq X^{\frac{1}{2}}$. Assume that $H$ is an integer which satisfies
\begin{align*}
1\leqslant H\ll\frac{X}{K}.
\end{align*}
From the arguments (4.16)-(4.24) in Tolev [18], we obtain
\begin{align*}
&\left|W(K)\right|^{2}
\ll\frac{X^{1+\varepsilon}}{H}\sum_{d_{1},d_2\leq D}\lambda(d_{1})\lambda(d_{2})\sum_{|h|<H}\left(1-\frac{|h|}{H}\right)\sum_{\frac{\mu X}{K_{1}}<\ell,\ell+h\leq\frac{X}{K}}\Lambda(\ell+h)\Lambda(\ell)\mathscr{F},
\end{align*}
where
\begin{align*}
\widetilde{K}=\max\bigg(K,\:\frac{\mu X}{\ell},\:\frac{\mu X}{\ell+h}\bigg),\quad\widetilde{K}_1=\min\bigg(K_1,\:\frac{X}{\ell},\:\frac{X}{\ell+h}\bigg),
\end{align*}
and
\begin{align*}
\mathscr{F}&:=\sum_{\widetilde{K}<k\leq\widetilde{K}_1\atop{d_1|k{\ell}+2,\,{d_{2}|k(\ell+h)+2}}}e\Bigl(k^c\bigl((\ell+h)^c-\ell^c\bigr)x\Bigr)\\
&=\sum_{\frac{\widetilde{K}-k_0}{[d_1,d_2]}<m\leq\frac{\widetilde{K}_1-k_0}{[d_1,d_2]}}e\Bigl(\bigl(k_0+m[d_1,d_2]\bigr)^c\bigl((\ell+h)^c-\ell^c\bigr)x\Bigr)\\
&=:\sum_{\frac{\widetilde{K}-k_0}{[d_1,d_2]}<m\leq\frac{\widetilde{K}_1-k_0}{[d_1,d_2]}}e(F(m)).
\end{align*}
For the case $h=0$ , we have $\mathscr{F}\ll K[d_{1},d_{2}]^{-1}.$ By the elementary estimate
\begin{align*}
\sum_{d_{1},d_{2}\leq D}[d_{1},d_{2}]^{-1}\ll {\mathcal{L}}^{3},
\end{align*}
we find that the contribution of $\mathscr{F}$ with $h=0$ to $|W(K)|^2$ is $\ll X^{2+\varepsilon}H^{-1}.$
\\
For the case $h\neq0,$ we have
\begin{align*}
\left|F^{(j)}(m)\right|\asymp|x|h\ell^{c-1}[d_1,d_2]K^{c-1}\cdot\left(K[d_1,d_2]^{-1}\right)^{1-j},\quad j\geqslant1.
\end{align*}
We apply Lemmas 2.4-2.5 with $(\kappa,\lambda)=(\frac{187}{659}+\varepsilon, \frac{374}{659}+\varepsilon)$ and find that
\begin{align*}
\mathscr{F}&\ll|x|^{-1}h^{-1}\ell^{1-c}[d_1,d_2]^{-1}K^{1-c}+\left(|x|h\ell^{c-1}[d_1,d_2]K^{c-1}\right)^{\frac{187}{659}+\varepsilon}\left(K[d_1,d_2]^{-1}\right)^{\frac{374}{659}+\varepsilon}\\
&\ll|x|^{-1}h^{-1}\ell^{1-c}[d_1,d_2]^{-1}K^{1-c}+|x|^{\frac{187}{659}}h^{\frac{187}{659}}\ell^{\frac{187}{659}(c-1)}[d_1,d_2]^{-\frac{187}{659}}K^{\frac{187c}{659}+\frac{187}{659}}.
\end{align*}
From the following estimate
\begin{align*}
\sum_{d_{1},d_{2}\leq D}[d_{1},d_{2}]^{-\frac{187}{659}}& \ll\sum_{d_{1},d_{2}\leq D}\left(\frac{(d_{1},d_{2})}{d_{1}d_{2}}\right)^{\frac{187}{659}}=\sum_{1\leq r\leq D}\sum_{k_{1}\leq\frac{D}{r}}\sum_{k_{2}\leq\frac{D}{r}}\frac{1}{r^{\frac{187}{659}}k_{1}^{\frac{187}{659}}k_{2}^{\frac{187}{659}}} \\
&\ll\sum_{1\leq r\leq D}r^{-{\frac{187}{659}}}\biggl(\sum_{k\leq\frac{D}{r}}k^{-{\frac{187}{659}}}\biggr)^{2}\ll\sum_{1\leq r\leq D}r^{-{\frac{187}{659}}}\biggl(\frac{D}{r}\biggr)^{\frac{944}{659}}\ll D^{{\frac{472}{659}}},
\end{align*}
we find the contribution of $\mathscr{F}$ with $h\neq0$ to $|W(K)|^2$ is
\begin{align*}
&\ll X^{1+\varepsilon}H^{-1}|x|^{-1}K^{1-c}\bigg(\sum_{d_{1},d_{2}\leq D}[d_{1},d_{2}]^{-1}\bigg)\bigg(\sum_{0<|h|<H}h^{-1}\bigg)\bigg(\sum_{\frac{\mu X}{K_{1}}<\ell,\ell+h\leq\frac{X}{K}}\ell^{1-c}\bigg) \\
&+X^{1+\varepsilon}H^{-1}|x|^{\frac{187}{659}}K^{\frac{187c}{659}+\frac{187}{659}}\bigg(\sum_{d_{1},d_{2}\leq D}[d_{1},d_{2}]^{-\frac{187}{659}}\bigg)\left(\sum_{0<|h|<H}h^{\frac{187}{659}}\right)\bigg(\sum_{\frac{\mu X}{K_1}<\ell,\ell+h\leq\frac{X}{K}}\ell^{\frac{187}{659}(c-1)}\bigg) \\
&\ll X^{3-c+\varepsilon}|x|^{-1}H^{-1}K^{-1}+X^{\frac{187c+1131}{659}}H^{\frac{187}{659}}D^{{\frac{472}{659}}}|x|^{\frac{187}{659}}K^{-\frac{285}{659}}.
\end{align*}
Then we get
\begin{align}
\left|W(K)\right|^2\ll X^{\varepsilon}\Big(X^2H^{-1}+X^{3-c}|x|^{-1}H^{-1}K^{-1}+X^{\frac{187c+1131}{659}}H^{\frac{187}{659}}D^{{\frac{472}{659}}}|x|^{\frac{187}{659}}K^{-\frac{285}{659}}\Big).
\end{align}
We choose
\begin{align*}
H_{0}=X^{\frac{187(1-c)}{846}}K^{\frac{285}{846}}D^{-\frac{472}{846}}|x|^{-\frac{187}{846}},\,\,\,\,\,\,\,\mbox{where}\,\, H=\Big[\min\{H_{0},XK^{-1}\}\Big].
\end{align*}
It is easy to see that
\begin{align}
H^{-1}\asymp H_{0}^{-1}+KX^{-1}.
\end{align}
By using (2.10) and (2.11), and noting that $K\geq X^{\frac{1}{2}}$ , we derive that
\begin{align*}
&\left|W(K)\right|^2\\
\ll\,\,& X^{\varepsilon}\Big(X^2\big(H_0^{-1}+KX^{-1}\big)+X^{3-c}|x\big|^{-1}K^{-1}\big(H_0^{-1}+KX^{-1}\big)+X^{\frac{187c+1131}{659}}H_0^{\frac{187}{659}}D^{{\frac{472}{659}}}|x|^{\frac{187}{659}}K^{-\frac{285}{659}}\Big)\\
\ll\,\,&  X^{\varepsilon}\Big(X^{\frac{187c}{846}+\frac{2725}{1692}}D^{\frac{472}{846}}|x|^{\frac{187}{846}}+X^{\frac{5}{3}}+X^{\frac{-659c}{846}+\frac{3571}{1692}}D^{\frac{472}{846}}|x|^{-\frac{659}{846}}+X^{2-c}|x|^{-1}\Big),
\end{align*}
which implies that
\begin{align*}
\left|W(K)\right|\ll X^{\varepsilon}\Big(X^{\frac{374c+2725}{3384}}D^{\frac{118}{423}}|x|^{\frac{187}{1692}}+X^{\frac{5}{6}}+X^{\frac{3571-1318c}{3384}}D^{\frac{118}{423}}|x|^{-\frac{659}{1692}}+X^{1-\frac{c}{2}}|x|^{-\frac{1}{2}}\Big).
\end{align*}
On the other hand, in the case $K<X^{\frac{1}{2}},$ we follow the argument (4.29)-(4.30) in Tolev [18] and derive that
\begin{align}
S_{2}^{''},\,S_3\ll X^{\varepsilon}\Big(X^{\frac{374c+2725}{3384}}D^{\frac{118}{423}}|x|^{\frac{187}{1692}}+X^{\frac{5}{6}}+X^{\frac{3571-1318c}{3384}}D^{\frac{118}{423}}|x|^{-\frac{659}{1692}}+X^{1-\frac{c}{2}}|x|^{-\frac{1}{2}}\Big).
\end{align}
Above all, from (2.5),(2.8)-(2.9) and (2.12), we deduce that
\begin{align*}
&L_1(x)\\
\ll \,&X^{\varepsilon}\Big(X^{\frac{11c}{82}+\frac{29}{41}}D^{\frac{25}{82}}|x|^{\frac{11}{82}}+X^{\frac{374c+2725}{3384}}D^{\frac{118}{423}}|x|^{\frac{187}{1692}}+X^{\frac{5}{6}}+X^{\frac{3571-1318c}{3384}}D^{\frac{118}{423}}|x|^{-\frac{659}{1692}}+X^{1-\frac{c}{2}}|x|^{-\frac{1}{2}}\Big)\\
\ll \,& X^{\varepsilon}\Big(X^{\frac{11c}{82}+\frac{29}{41}+\frac{25}{82}\delta}+X^{\frac{374c+2725}{3384}+\frac{118}{423}\delta}+X^{\frac{5}{6}}+X^{\frac{3571}{3384}+\frac{118}{423}\delta-\frac{659}{1692}\xi}+X^{1-\frac{1}{2}\xi}\Big),
\end{align*}
and thus, for $1<c<\frac{1787}{1502}$, there holds
\begin{align*}
\sup_{\tau\leq|x|\leq K}\left|L_{1}(x)\right|\ll X^{\frac43-\frac13c-\varepsilon},
\end{align*}
which completes the proof of Lemma 2.10.\\

{\bf Lemma 2.11}.
\mbox{(i)} {\it Define
\begin{align*}
T(x)=\sum_{d\leq D}\sum_{\mu X<n\leq X\atop{d|n+2}}e(n^cx).
\end{align*}
Then for $0<|x|\leq 2K,$ we have}
\begin{align*}
T(x)\ll X^{\frac{13c}{31}+\frac{3}{31}+\varepsilon}D^{\frac{13}{31}}+|x|^{-1}X^{1-c}{\mathcal L}.
\end{align*}
\mbox{(ii)} {\it For any Borel measurable bounded function $G$ on $[\tau,K],$ we have}
\begin{align*}
&\int_{\tau\leq|x|\leq K}L(x)G(x)\mathrm{d}x\\
\ll \,&X^{\frac{13c}{62}+\frac{17}{31}+\frac{13}{62}\delta+\varepsilon}\times\left(\intop_{\tau\leq|x|\leq K}|G(x)|\mathrm{d}x\right)+X^{1-\frac{c}{2}+\frac{\varepsilon}{2}}\left(\max_{\tau\leq|x|\leq K}|G(x)|\right)^{\frac12}\left(\intop_{\tau\leq|x|\leq K}|G(x)|\mathrm{d}x\right)^{\frac12}.
\end{align*}
\mbox{(iii)} {\it We have}
\begin{align*}
&\int_{\tau\leq|x|\leq K}|L(x)|^4\mathrm{d}x \ll X^{4-c-\varepsilon}.
\end{align*}

{\bf Proof}. For (i), it is obvious that
\begin{align}
T(x)&=\sum_{d\leq D}\sum_{\frac{\mu X+2}{d}<h\leq\frac{X+2}{d}}e\big((hd-2)^{c}x\big)\notag\\
&\ll \mathcal{L}\max_{\mathscr{D}\leq D}\sum_{D\sim \mathscr{D}}\bigg|\sum_{\frac{\mu X+2}{d}<h\leq\frac{X+2}{d}}e\big((hd-2)^{c}x\big)\bigg|.
\end{align}
For the inner sum in (2.13), it follows from Lemma 2.4 with the exponential pair $(\kappa,\lambda)=BA^4B(0,1)=(\frac{13}{31},\frac{16}{31})$ that
\begin{align}
\sum_{\frac{\mu X+2}d<h\leq\frac{X+2}d}e\big((hd-2)^cx\big)&\ll\big(|x|dX^{c-1}\big)^{\frac{13}{31}}\bigg(\frac Xd\bigg)^{\frac{16}{31}}+\frac{X^{1-c}}{|x|d}\notag\\
&\ll|x|^{\frac{13}{31}}d^{-\frac{3}{31}}X^{\frac {13c}{31}+\frac{3}{31}}+\frac{X^{1-c}}{|x|d}.
\end{align}
On inserting (2.14) into (2.13), we obtain
\begin{align*}
T(x)& \ll{\mathcal L}\max_{\mathscr{D}\leq D}\sum_{D\sim \mathscr{D}}\left(|x|^{\frac{13}{31}}d^{-\frac{3}{31}}X^{\frac {13c}{31}+\frac{3}{31}}\right) \ll X^{\frac{13c}{31}+\frac{3}{31}+\varepsilon}D^{\frac{13}{31}}+|x|^{-1}X^{1-c}{\mathcal L},
\end{align*}
which completes the proof of (i). Next we turn to the proof of (ii).
\begin{align*}
&\left|\intop_{\tau\leq|x|\leq K}L(x)G(x)\mathrm{d}x\right|
=\,\,\left|\sum_{d\leq D}\lambda(d)\sum_{\mu X<p\leq X\atop{d|p+2}}(\log p)\intop_{\tau\leq|x|\leq K}e(p^{c}x)G(x)\mathrm{d}x\right| \\
\leq\,\,&\sum_{d\leq D}\sum_{\mu X<p\leq X\atop{d|p+2}}(\log p)\left|\intop_{\tau\leq|x|\leq  K}e(p^{c}x)G(x)\mathrm{d}x\right| \,\,\leq \mathcal{L}\sum_{d\leq D}\sum_{\mu X<n\leq X\atop{d|n+2}}\left|\intop_{\tau\leq|x|\leq  K}e(n^{c}x)G(x)\mathrm{d}x\right|.
\end{align*}
By Cauchy's inequality, we have
\begin{align}
&\left|\intop_{\tau\leq|x|\leq  K}L(x)G(x)\mathrm{d}x\right|^{2} \notag\\
\ll \,\,&X\mathcal{L}^{2}\sum_{\mu X<n\leq X}\left(\sum_{d\leq D\atop{d|n+2}}\left|\intop_{\tau\leq|x|\leq  K}e(n^{c}x)G(x)\mathrm{d}x\right|\right)^{2} \notag\\
\ll \,\,&X\mathcal{L}^{2}\sum_{\mu X<n\leq X}\bigg(\sum_{d\leq D\atop{d|n+2}}1\bigg)\sum_{d\leq D\atop{d|n+2}}\left|\intop_{\tau\leq|x|\leq K}e(n^{c}x)G(x)\mathrm{d}x\right|^{2} \notag\\
\ll \,\,&X^{1+\varepsilon}\sum_{d\leq D}\sum_{\mu X<n\leq X\atop d|n+2}\left|\intop_{\tau\leq|x|\leq  K}e(n^{c}x)G(x)\mathrm{d}x\right|^{2} \notag\\
=\,\,&X^{1+\varepsilon}\sum_{d\leq D}\sum_{\mu X<n\leq X\atop{d|n+2}}\left(\intop_{\tau\leq|x|\leq K}e(n^{c}x)G(x)\mathrm{d}x\right)\left(\intop_{\tau\leq|y|\leq  K}\overline{e(n^{c}y)G(y)}\mathrm{d}y\right) \notag\\
\ll \,\,&X^{1+\varepsilon}\intop_{\tau\leq|y|\leq  K}\left|G(y)\right|\mathrm{d}y\intop_{\tau\leq|x|\leq  K}\left|G(x)\right|\left|T(x-y)\right|\mathrm{d}x.\notag
\end{align}
On combining Lemma 2.11 (i) and the trivial estimate $|T(x-y)|\ll X \mathcal{L}$, we obtain
\begin{align}
&\left|\intop_{\tau\leq|x|\leq  K}L(x)G(x)\mathrm{d}x\right|^{2} \notag\\
\ll \,\,&X^{\frac{13c}{31}+\frac{34}{31}+\frac{13}{31}\delta+2\varepsilon}\intop_{\tau\leq|y|\leq K}\left|G(y)\right|\mathrm{d}y\intop_{\tau\leq|x|\leq K}\left|G(x)\right|\mathrm{d}x \notag\\
\,+\,&
X^{1+\varepsilon}\mathcal{L}\intop_{\tau\leq|y|\leq  K}\left|G(y)\right|\mathrm{d}y\intop_{\tau\leq|x|\leq  K}\left|G(x)\right|\min\left(X, \frac{X^{1-c}}{|x-y|}\right)\mathrm{d}x.
\end{align}
It is easy to see that
\begin{align}
&\intop_{\tau\leq|x|\leq K}|G(y)|\intop_{\tau\leq|x|\leq  K}|G(x)|\min\left(X,\frac{X^{1-c}}{|x-y|}\right)\mathrm{d}x\mathrm{d}y \notag\\
\leq\,\,&\frac{1}{2}\int\limits_{\tau\leq|y|\leq  K}\int\limits_{\tau\leq|x|\leq K}(|G(x)|^{2}+|G(y)|^{2})\min\left(X,\frac{X^{1-c}}{|x-y|}\right)\mathrm{d}x\mathrm{d}y \notag\\
=\,\,&\intop_{\tau\leq|x|\leq  K}|G(x)|^2\intop_{\tau\leq|y|\leq K}\min\left(X,\frac{X^{1-c}}{|x-y|}\right)\mathrm{d}y\mathrm{d}x\notag\\
=\,\,&\int\limits_{\tau\leq|x|\leq  K}|G(x)|^2\left(X\int\limits_{\tau\leq|y|\leq K\atop{|x-y|\leq X^{-c}}}dy+\int\limits_{\tau\leq|y|\leq K\atop{X^{-c}\leq|x-y|\leq2K}}\frac{X^{1-c}}{|x-y|}\mathrm{d}y\right)\mathrm{d}x \notag\\
\ll\,\,& X^{1-c}\mathcal{L}\intop_{\tau\leq|x|\leq K}|G(x)|^2dx.
\end{align}
On combining (2.15) and (2.16), we have
\begin{align*}
&\left|\intop_{\tau\leq|x|\leq K}L(x)G(x)\mathrm{d}x\right|^{2} \notag\\
\ll \,\,&X^{\frac{13c}{31}+\frac{34}{31}+\frac{13}{31}\delta+2\varepsilon}\left(\intop_{\tau\leq|x|\leq  K}\left|G(x)\right|\mathrm{d}x\right)^2 + X^{2-c+\varepsilon}\max\limits_{\tau\leq|y|\leq K}\left|G(x)\right|\times\intop_{\tau\leq|x|\leq  K}\left|G(x)\right|\mathrm{d}x.
\end{align*}
Hence the proof of Lemma 2.11 (ii) is proved.

For (iii), it follows from Lemmas 2.8 (v), 2.10 and 2.11 (ii) with $G(x)=\overline{L(x)}|L(x)|^2$ that
\begin{align*}
\int_{\tau\leq|x|\leq K}|L(x)|^4\mathrm{d}x &\ll X^{\frac{13c}{62}+\frac{17}{31}+\frac{13}{62}\delta+\varepsilon}X^{2}+X^{1-\frac{c}{2}+\frac{\varepsilon}{2}}(X^{\frac43-\frac13c-\varepsilon})^{\frac32}X\\
&\ll X^{4-c-\varepsilon},
\end{align*}
which now completes the proof of Lemma 2.11.\\

{\bf Lemma 2.12}. {\it We have}
\begin{align*}
\mathbb{A}=\max_{R_{2}\in(N,2N]}\int_{N}^{2N}\left|\int_{\tau \le|x|\leq K}e\left(x(R_{1}-R_{2})\right)\mathrm{d}x\right|\mathrm{d}R_{1}\ll\log N.
\end{align*}

{\bf Proof}. See  Laporta [10, Lemma 1].\\

{\bf Lemma 2.13}. {\it Let $\Omega_1$ and $\Omega_2$ be measurable subsets of $\mathbb{R}^{n}$ and let ${\left\|f\right\|_{j}=}$
$\left(\int_{\Omega_j}|f(y)|^2\mathrm{d}y\right)^{\frac{1}{2}},$ $\langle f,g\rangle_j=\int_{\Omega_j}f(y)\overline{g(y)}\mathrm{d}y$ denote the usual norm and inner product in $L^{2}(\Omega_{j},\mathbb{C})\,(j=1,2)$ respectively. Let $\zeta\in L^{2}(\Omega_{1},\mathbb{C}),\xi\in L^{2}(\Omega_{2},\mathbb{C}),$ and let $\omega$ be a measurable complex-valued function on $\Omega_{1}\times\Omega_{2}$ such that}
\begin{align*}
\sup_{x\in\Omega_{1}}\int_{\Omega_{2}}|\omega(x,y)|\mathrm{d}y<+\infty,\,\,\,\sup_{y\in\Omega_{2}}\int_{\Omega_{1}}|\omega(x,y)|\mathrm{d}x<+\infty.
\end{align*}
Then we have
\begin{align*}
\left|\int_{\Omega_1}\zeta(x)\langle\xi,\omega(x,\cdot)\rangle_2\mathrm{d}x\right|\ll\|\xi\|_2\|\zeta\|_1\left(\sup_{x^{'}\in\Omega_1}\int_{\Omega_1}\left|(\omega(x,\cdot),\omega(x^{'},\cdot))_2\right|\mathrm{d}x\right)^{\frac{1}{2}}.
\end{align*}

{\bf Proof}. See  Laporta [10, Lemma 2].

\section{Proof of the Theorem 1}
\setcounter{equation}{0}
Let $N\le R \le 2N$ and
\begin{align}
&\Gamma(R)=\sum_{\mu X<p_1,p_2\leqslant X\atop{|p_1^c+p_2^c-R|<\Delta\atop{(p_{i}+2,P(z))=1,\,i=1,2}}}(\log p_{1})(\log p_{2})\notag
\end{align}
By the definition of $\varphi(y)$ (see Lemma 2.1) and the inverse Fourier transformation formula, we have
\begin{align*}
&\Gamma(R) \ge \sum_{\mu X<p_1,p_2\leq X\atop{(p_i+2,P(z))=1,\,i=1,2}}(\log p_{1})(\log p_{2})\varphi(p_1^c+p_2^c-R)\notag\\
&\,\,\,\,\,\,\,\,\,\,\,\,\,\,=\,\sum_{\mu X<p_{1},p_{2}\leq X}(\log p_{1})(\log p_{2})\Lambda_{1}\Lambda_{2}\varphi(p_{1}^{c}+p_{2}^{c}-R)\\
&\,\,\,\,\,\,\,\,\,\,\,\,\,\,=: \widetilde{\Gamma}(R)
\end{align*}
From the linear sieve and Lemma 2.3 (i), we obtain
\begin{align*}
\widetilde{\Gamma}(R)
&\geq\sum_{\mu X<p_{1},p_{2}\leq X}(\log p_{1})(\log p_{2})\varphi(p_{1}^{c}+p_{2}^{c}-R)(\Lambda_{1}^{-}\Lambda_{2}^{+}+\Lambda_{1}^{+}\Lambda_{2}^{-} -\Lambda_{1}^{+}\Lambda_{2}^{+})\\
&=: \Gamma_{1}(R)+\Gamma_{2}(R)- \Gamma_{3}(R),
\end{align*}
say. By the symmetric property, it is clear that
\begin{align}
&\Gamma_{1}(R)=\Gamma_{2}(R)
=\sum_{\mu X<p_{1},p_{2}\leq X}(\log p_{1})(\log p_{2})
\Lambda_{1}^{-}\Lambda_{2}^{+}\varphi\left(p_{1}^{c}+p_{2}^{c}-R\right), \\
&\Gamma_{3}(R)=\sum_{\mu X<p_{1},p_{2}\leq X}(\log p_{1})(\log p_{2})\Lambda_{1}^{+}\Lambda_{2}^{+}\varphi\left(p_{1}^{c}+p_{2}^{c}-R\right),
\end{align}
and thus
\begin{align}
\Gamma(R)\ge2\Gamma_1(R)-\Gamma_3(R).
\end{align}

We first consider the sum $\Gamma_{1}(R)$. According to the Fourier's inverse transformation, we have
\begin{align}
\Gamma_1(R)
&=\,\sum_{\mu X<p_{1},p_{2}\leq X}(\log p_{1})(\log p_{2})\Lambda_{1}^{-}\Lambda_{2}^{+}\int_{-\infty}^{+\infty}\Theta(x)e\left((p_{1}^{c}+p_{2}^{c}-R)x\right)\mathrm{d}x \notag\\
&=\,\int_{-\infty}^{+\infty}L^{-}(x)L^{+}(x)\Theta(x)e(-Rx)\mathrm{d}x \notag\\
&=\,\Bigg(\,\intop_{|x|> K}+\intop_{\tau \le |x| \le K}+\intop_{|x|>\tau}\Bigg)L^{-}(x)L^{+}(x)\Theta(x)e(-Rx)\mathrm{d}x \notag\\
&=:\,\Gamma_{1}^{(1)}(R)+\Gamma_{1}^{(2)}(R)+\Gamma_{1}^{(3)}(R)\notag.
\end{align}
Similarly, for the quantity $\Gamma_3(R)$, we also have
\begin{align}
\Gamma_3(R)& =\int_{-\infty}^{+\infty}(L^{+}(x))^{2}\Theta(x)e(-Rx)\mathrm{d}x  \notag\\
&=\Bigg(\,\intop_{|x|> K}+\intop_{\tau \le |x| \le K}+\intop_{|x|>\tau}\Bigg)(L^{+}(x))^{2}\Theta(x)e(-Rx)\mathrm{d}x \notag\\
&=:\Gamma_{3}^{(1)}(R)+\Gamma_{3}^{(2)}(R)+\Gamma_{3}^{(3)}(R)\notag.
\end{align}
Define
\begin{align*}
H(R)&=\int_{-\infty}^{+\infty}I^2(x)\Theta(x)e(-Rx)\mathrm{d}x,\notag\\
H_1(R)&=\int_{-\tau}^{\tau}I^{2}(x)\Theta(x)e(-Rx)\mathrm{d}x.
\end{align*}
Next, we give the estimation over the six terms in $\Gamma_1(R)$ and $\Gamma_3(R)$ respectively by giving the following Propositions.\\
{\bf Proposition 3.1}. {\it For $1<c<\frac{1787}{1502},$ we have}
\begin{align*}
\int_{N}^{2N}|\Gamma_{1}^{(1)}(R)|^{2}\mathrm{d}R\ll N,\,\,\,\,\,\,\,\,\,\,\,\,\,\,\,\,\,\,\int_{N}^{2N}|\Gamma_{3}^{(1)}(R)|^{2}\mathrm{d}R\ll N.
\end{align*}
{\bf Proof}.
By the trivial estimate $L^\pm(x)\ll X^{1+\varepsilon}$ and Lemma 2.1, we get
\begin{align*}
&\int_{N}^{2N}|\Gamma_{1}^{(1)}(R)|^{2}\mathrm{d}R \ll N\left|\int_{K}^{+\infty}|L^\pm(x)|^{2}|\Theta(x)|\mathrm{d}x\right|^{2} \\
\ll \,\,&NX^{4+\varepsilon}\left|\int_{K}^{+\infty}\left(\frac{r}{2\pi x\Delta}\right)^{r}\frac{\mathrm{d}x}{\pi x}\right|^{2}\ll NX^{4+\varepsilon}\left(\frac{r}{2\pi K\Delta}\right)^{2r}\ll N.
\end{align*}
\\
{\bf Proposition 3.2}. {\it For $1<c<\frac{1787}{1502},$ we have}
\begin{align*}
&\int_{N}^{2N}|\Gamma_{1}^{(2)}(R)|^{2}\mathrm{d}R\ll \Delta^2\mathcal{L} N^{\frac{4}{c}-1-\varepsilon},\\
&\int_{N}^{2N}|\Gamma_{3}^{(2)}(R)|^{2}\mathrm{d}R\ll \Delta^2\mathcal{L} N^{\frac{4}{c}-1-\varepsilon}.
\end{align*}
{\bf Proof}. For $\Gamma_{1}^{(2)}(R)$,  we set
\begin{align*}
&\Omega_{1}=\{R:R\sim N\},\,\,\,\,\,\,\,\Omega_{2}=\{x:\tau \le|x|\le K\},\,\,\,\,\,\,\,\xi(x)=L^{-}(x)L^{+}(x)\Theta(x),\\
&\omega(x,R)=e(Rx),\,\,\,\,\,\,\,\,\,\,\,\,\zeta(R)=\overline{\Gamma_{1}^{(2)}(R)}.
\end{align*}
Then it follows from Lemma 2.13 that
\begin{align}
&\int_{N}^{2N}|\Gamma_{1}^{(2)}(R)|^{2}\mathrm{d}R \notag\\
=\,\,&\int_{\Omega_{1}}\overline{\Gamma_{1}^{(2)}(R)}\langle L^{-}(x)L^{+}(x)\Theta(x),e(Rx)\rangle_{2}\mathrm{d}R \notag\\
\ll\,\,&\left(\int_{\Omega_2}|{L^{-}(x)}^2\Theta(x)|^2\mathrm{d}x\right)^{\frac{1}{4}}\left(\int_{\Omega_2}|{L^{+}(x)}^2\Theta(x)|^2\mathrm{d}x\right)^{\frac{1}{4}}\left(\int_{\Omega_1}\Big|\overline{\Gamma_{1}^{(2)}(R)}\Big|^2\mathrm{d}R\right)^{\frac{1}{2}} \notag\\
\times&\left(\sup_{x^{'}\in\Omega_1}\int_{\Omega_1}\left|\langle e(Rx),e(Rx^{'})\rangle_2\right|\mathrm{d}x\right)^{\frac{1}{2}}.
\end{align}
The inequality (3.4) together with Lemmas 2.1, 2.11 (iii), 2.12 and  Cauchy-Schwarz inequality yields that
\begin{align*}
\int_N^{2N}|\Gamma_{1}^{(2)}(R)|^{2}\mathrm{d}R
&\ll\mathbb{A}\left(\intop_{\tau\leq|x|\leq K}|L^{-}(x)|^4|\Theta(x)|^2\mathrm{d}x\right)^{\frac 12}\left(\intop_{\tau\leq|x|\leq K}|L^{+}(x)|^4|\Theta(x)|^2\mathrm{d}x\right)^{\frac 12}\\
&\ll\Delta^2\mathcal{L}\left(\intop_{\tau\leq|x|\leq K}\big|L^{-}(x)\big|^{4}\mathrm{d}x\right)^{\frac{1}{2}}\left(\intop_{\tau\leq|x|\leq K}\big|L^{+}(x)\big|^{4}\mathrm{d}x\right)^{\frac{1}{2}}\\
&\ll\Delta^2\mathcal{L} X^{4-c-\varepsilon}\,\,\ll\Delta^2\mathcal{L} N^{\frac{4}{c}-1-\varepsilon}.
\end{align*}
Similarly, we have
\begin{align*}
\int_N^{2N}|\Gamma_{3}^{(2)}(R)|^{2}\mathrm{d}R \ll \Delta^2\mathcal{L}\intop_{\tau\leq|x|\leq  K}\left|L^{+}(x)\right|^{4}\mathrm{d}x\ll \Delta^2\mathcal{L} N^{\frac{4}{c}-1-\varepsilon}.
\end{align*}
\\
{\bf Proposition 3.3}. {\it For $1<c<\frac{1787}{1502},$ we have}
\begin{align*}
&\int_{N}^{2N}\left|\Gamma_{1}^{(3)}(R) - \mathcal{M}^{-}\mathcal{M}^{+}H_1(R)\right|^2\mathrm{d}R\ll \Delta^2{\mathcal{L}}^{9-2A}N^{\frac{4}{c}-1},\\
&\int_{N}^{2N}\left|\Gamma_{3}^{(3)}(R) - \mathcal{M}^{-}\mathcal{M}^{+}H_1(R)\right|^2\mathrm{d}R\ll \Delta^2{\mathcal{L}}^{9-2A}N^{\frac{4}{c}-1}.
\end{align*}
{\bf Proof}. We give the evaluation of the first inequality, while the second one is similar.
It is easy to see
\begin{align}
&\left|\Gamma_{1}^{(3)}(R) - \mathcal{M}^{-}\mathcal{M}^{+}H_1(R)\right|^2 \notag\\
=\,\,& \left|\int_{-\tau}^{\tau}\left(L^{-}(x)L^{+}(x)-\mathcal{M}^{-}\mathcal{M}^{+}I^{2}(x)\right)\Theta(x)e(-Nx)\mathrm{d}x\right|^2 \notag\\
=\,\,& \left|\int_{-\tau}^{\tau}\left(L^{+}(x)\big(L^{-}(x)-\mathcal{M}^{-}I(x)\big)+I(x)\mathcal{M}^{-}\big(L^{+}(x)-\mathcal{M}^{+}I(x)\big)\right)\Theta(x)e(-Nx)\mathrm{d}x\right|^2\\
\ll\,\,&\Delta^2 X^2{\mathcal{L}}^{2-2A}\int_{-\tau}^{\tau}\Big(\overline{L^{+}(x)}+\overline{I(x)}\Big)\mathrm{d}x
\int_{-\tau}^{\tau}\Big(L^{+}(y)+I(y)\Big)e((x-y)R)\mathrm{d}y,
\end{align}
where (2.4), Lemmas 2.1 and  2.7 (ii) are used in (3.5)-(3.6). Then it follows from Lemma 2.8 (i) and (iii) that
\begin{align*}
&\int_{N}^{2N}\left|\Gamma_{1}^{(3)}(R) - \mathcal{M}^{-}\mathcal{M}^{+}H_1(R)\right|^2\mathrm{d}R \notag\\
\ll\,\,&\Delta^2 X^2{\mathcal{L}}^{2-2A}\int_{-\tau}^{\tau}\Big(\overline{L^{+}(x)}+\overline{I(x)}\Big)\mathrm{d}x
\int_{-\tau}^{\tau}\Big(L^{+}(y)+I(y)\Big)\left(\int_{N}^{2N}e((x-y)R)\mathrm{d}R\right)\mathrm{d}y\\
\leq\,\,&\frac{1}{2}\Delta^2 X^2{\mathcal{L}}^{2-2A}\int_{-\tau}^{\tau}\int_{-\tau}^{\tau}\left(\Big(\overline{L^{+}(x)}+\overline{I(x)}\Big)^{2}+\Big(L^{+}(y)+I(y)\Big)^{2}\right)\min\left(N,\frac{1}{|x-y|}\right)\mathrm{d}y\mathrm{d}x \notag\\
\ll\,\,&\Delta^2 X^2{\mathcal{L}}^{2-2A}\int_{-\tau}^{\tau}\left(|L(x)|^2 + |I(x)|^2\right)\int_{-\tau}^{\tau}\min\left(N,\frac{1}{|x-y|}\right)\mathrm{d}y\mathrm{d}x\notag\\
\ll\,\,&\Delta^2 X^2{\mathcal{L}}^{2-2A}\int_{-\tau}^{\tau}\left(|L(x)|^2 + |I(x)|^2\right)\Biggr(N\int\limits_{|y|\leq \tau\atop{|x-y|\leq N^{-1}}}dy+\int\limits_{|y|\leq \tau\atop{N^{-1}\leq|x-y|\leq2\tau}}\frac{1}{|x-y|}\mathrm{d}y\Biggr)\mathrm{d}x \notag\\
\ll\,\,&\Delta^2 X^{4-c}{\mathcal{L}}^{9-2A}\,\,\ll \Delta^2{\mathcal{L}}^{9-2A}N^{\frac{4}{c}-1}.
\end{align*}
\\
{\bf Proposition 3.4}. {\it For $1<c<\frac{1787}{1502},$ we have}
\begin{align*}
&\int_{N}^{2N}\left|\mathcal{M}^{-}\mathcal{M}^{+}H_1(R) - \mathcal{M}^{-}\mathcal{M}^{+}H(R)\right|^2\mathrm{d}R\ll \Delta^2{\mathcal{L}}^{4}N^{\frac{5}{c}-2}.
\end{align*}
{\bf Proof}. Since $\left|\frac{\mathrm{d}}{\mathrm{d}t}(t^cx)\right|\gg|x|X^{c-1}$ for $t\in(\mu X,X]$, then on applying Titchmarsh [17, Lemma 4.2], we find that $I(x) \ll |x|^{-1}X^{1-c}$. Then it follows from (2.4) and Lemma 2.1 that
\begin{align*}
&\int_{N}^{2N}\left|\mathcal{M}^{-}\mathcal{M}^{+}H_1(R) - \mathcal{M}^{-}\mathcal{M}^{+}H(R)\right|^2\mathrm{d}R\notag\\
\ll\,\,&{\mathcal{L}}^4\int_{N}^{2N}\int_\tau^{+\infty}\left|I(x)\right|^4\big|\Theta(x)\big|^2\mathrm{d}x\mathrm{d}R\notag\\
\ll\,\,&\Delta^2N{\mathcal{L}}^4 X^{4-4c}\int_\tau^{+\infty}\frac{1}{x^4}\mathrm{d}x
\,\,\ll\,\Delta^2N{\mathcal{L}}^4 X^{4-4c}\tau^{-3}\notag\\
\ll\,\,&\Delta^2N{\mathcal{L}}^4 X^{4-c-3\xi}\,\,\ll \Delta^2{\mathcal{L}}^{4}N^{\frac{5}{c}-2}.
\end{align*}
Combining Propositions 3.1-3.4, we have
\begin{align}
&\int_{N}^{2N}|\Gamma_{1}(R)-\mathcal{M}^{-}\mathcal{M}^{+}H(R)|^{2}\mathrm{d}R \notag\\
\leq\,\,&\int_{N}^{2N}|\Gamma_{1}^{(1)}(R)|^{2}\mathrm{d}R+\int_{N}^{2N}|\Gamma_{1}^{(2)}(R)|^{2}\mathrm{d}R \notag\\
+&\int_{N}^{2N}\left|\Gamma_{1}^{(3)}(R) - \mathcal{M}^{-}\mathcal{M}^{+}H_1(R)\right|^2\mathrm{d}R+\int_{N}^{2N}\left|\mathcal{M}^{-}\mathcal{M}^{+}H_1(R) - \mathcal{M}^{-}\mathcal{M}^{+}H(R)\right|^2\mathrm{d}R\notag\\
\ll\,\,&\Delta^2{\mathcal{L}}^{9-2A}N^{\frac{4}{c}-1}.
\end{align}
By (3.7), we find that for $1<c<\frac{1787}{1502}$, for all $R \in (N,2N]\setminus \mathcal{A}$ with $|\mathcal{A}|=O(N{\mathcal{L}}^{\frac92-A})$, we have
\begin{align*}
\Gamma_{1}(R) = \mathcal{M}^{-}\mathcal{M}^{+}H(R) + O\left(\Delta{\mathcal{L}}^{\frac92-A}N^{\frac{2}{c}-1}\right).
\end{align*}
On recalling (2.4), (3.1)-(3.3), Lemmas 2.2 (ii) and 2.9 (i), we obtain
\begin{align*}
\Gamma(R) &\geq \big(2\mathcal{M}^{-}-\mathcal{M}^{+}\big)\mathcal{M}^{+}H(R)+ O\left(\Delta{\mathcal{L}}^{\frac92-A}N^{\frac{2}{c}-1}\right)\notag\\
&\geq \left(2f\bigg(\frac{\log D}{\log z_1}\bigg)-F\bigg(\frac{\log D}{\log z_1}\bigg)\bigg)\big(1+O({\mathcal L}^{-1/3})\big)\mathfrak{P}^2H(R)+ O\left(\Delta{\mathcal{L}}^{\frac92-A}N^{\frac{2}{c}-1}\right)\right.  \notag\\
&=\left(2f\bigg(\frac{53}{20}\bigg)-F\bigg(\frac{53}{20}\bigg)\bigg)\mathfrak{P}^2H(R)+O\big(\Delta{\mathcal{L}}^{-\frac73}R^{\frac{2}{c}-1}\big)\right.  \notag\\
&=\frac{80e^\gamma}{53}\bigg(\log\frac{33}{20}-\frac12\bigg)\mathfrak{P}^2H(R)+O\big(\Delta{\mathcal{L}}^{-\frac73}R^{\frac{2}{c}-1}\big) \notag\\
&\gg \Delta{\mathcal{L}}^{-2}R^{\frac{2}{c}-1}.
\end{align*}
Therefore, we find that for $1<c<\frac{1787}{1502}$, for all $R \in (N,2N]\setminus \mathcal{A}$ with $|\mathcal{A}|=O\bigg(N\left(\frac{\log N}{c}\right)^{\frac92-A}\bigg)$, $\Gamma(R)>0$ and the inequality (1.3) has a solution in primes $p_1,p_2$ such that each of the numbers $p_{i}+2\,(i=1,2)$ had at most [$\eta_1^{-1}]=\left[\frac{79606}{35740-30040c}\right]$  prime factors. This completes the proof of the Theorem 1.

\section{Proof of the Theorem 2}
\setcounter{equation}{0}
Let
\begin{align}
&\mathcal{B}(N)=\sum_{\mu X<p_1,p_2,p_3,p_4\leqslant X\atop{|p_1^c+p_2^c+p_3^c+p_4^c-N|<\Delta\atop{(p_{i}+2,P(z))=1,\,i=1,2,3,4}}}\left(\prod_{j=1}^{4}\log p_{j}\right)\notag\\
&\mathcal{B}^{'}(N)
=\,\sum_{\mu X<p_1,p_2,p_3,p_4\leq X\atop{(p_i+2,P(z))=1,\,i=1,2,3,4}}\left(\prod_{j=1}^{4}\log p_{j}\right)\varphi\:(p_1^c+p_2^c+p_3^c+p_4^c-N)\notag\\
&\,\,\,\,\,\,\,\,\,\,\,\,\,\,\,\,=\,\sum_{\mu X<p_{1},p_{2},p_{3},p_{4}\leq X}\left(\prod_{j=1}^{4}\log p_{j}\right)\Lambda_{1}\Lambda_{2}\Lambda_{3}\Lambda_{4}\varphi(p_{1}^{c}+p_{2}^{c}+p_{3}^{c}+p_{4}^{c}-N).
\end{align}
By the definition of $\varphi(y)$\,(see Lemma 2.1), we have
\begin{align*}
\mathcal{B}(N)\geq\mathcal{B}^{'}(N).
\end{align*}
From the linear sieve, we have $\Lambda_{i}^{-}\leq\Lambda_{i}\leq\Lambda_{i}^{+}$. Furthermore, on substituting the quantity from the right-hand side of  Lemma 2.3 (ii) for $\Lambda_{1}\Lambda_{2}\Lambda_{3}\Lambda_{4}$ in (4.1), we find that
\begin{align*}
\mathcal{B}^{'}(N) \geq \mathcal{B}_1(N)+\mathcal{B}_2(N)+\mathcal{B}_3(N)+\mathcal{B}_4(N)-3\mathcal{B}_5(N),
\end{align*}
where $\mathcal{B}_1(N),...,\mathcal{B}_5(N)$ are the contributions coming from the consecutive terms of the right-hand side of  Lemma 2.3 (ii). It is clear that
\begin{align*}
&\mathcal{B}_1(N)=\mathcal{B}_2(N)=\mathcal{B}_3(N)=\mathcal{B}_4(N)\notag\\
=\,\,&\sum_{\mu X<p_{1},p_{2},p_{3},p_{4}\leq X}\left(\prod_{j=1}^{4}\log p_{j}\right)
\Lambda_{1}^{-}\Lambda_{2}^{+}\Lambda_{3}^{+}\Lambda_{4}^{+}\varphi\left(p_{1}^{c}+p_{2}^{c}+p_{3}^{c}+p_{4}^{c}-N\right), \\
&\mathcal{B}_5(N)=\,\sum_{\mu X<p_{1},p_{2},p_{3},p_{4}\leq X}\left(\prod_{j=1}^{4}\log p_{j}\right)\Lambda_{1}^{+}\Lambda_{2}^{+}\Lambda_{3}^{+}\Lambda_{4}^{+}\varphi\left(p_{1}^{c}+p_{2}^{c}+p_{3}^{c}+p_{4}^{c}-N\right),
\end{align*}
and thus
\begin{align}
\mathcal{B}(N)\ge4\mathcal{B}_1(N)-3\mathcal{B}_5(N).
\end{align}

We first consider the sum $\mathcal{B}_1(N)$. According to the Fourier's inverse transformation, we have
\begin{align}
\mathcal{B}_1(N)
=\,\,&\sum_{\mu X<p_{1},p_{2},p_{3},p_{4}\leq X}\left(\prod_{j=1}^{4}\log p_{j}\right)\Lambda_{1}^{-}\Lambda_{2}^{+}\Lambda_{3}^{+}\Lambda_{4}^{+}\int_{-\infty}^{+\infty}\Theta(x)e\Big((p_{1}^{c}+p_{2}^{c}+p_{3}^{c}+p_{4}^{c}-N)x\Big)\mathrm{d}x \notag\\
=\,\,&\int_{-\infty}^{+\infty}L^{-}(x)(L^{+}(x))^{3}\Theta(x)e(-Nx)\mathrm{d}x.
\end{align}

We divide the integral (4.3) into three parts as follows:
\begin{align}
\mathcal{B}_1(N)=\mathcal{B}_{1}^{(1)}(N)+\mathcal{B}_{1}^{(2)}(N)+\mathcal{B}_{1}^{(3)}(N),
\end{align}
where
\begin{align*}
\mathcal{B}_{1}^{(1)}(N) &= \intop_{|x|> K}L^{-}(x)(L^{+}(x))^{3}\Theta(x)e(-Nx)\mathrm{d}x,\\
\mathcal{B}_{1}^{(2)}(N) &= \intop_{\tau \le|x|\le K}L^{-}(x)(L^{+}(x))^{3}\Theta(x)e(-Nx)\mathrm{d}x,\\
\mathcal{B}_{1}^{(3)}(N) &= \intop_{|x|<\tau}L^{-}(x)(L^{+}(x))^{3}\Theta(x)e(-Nx)\mathrm{d}x.
\end{align*}
Similarly, for the quantity $\mathcal{B}_5(N)$, we also have
\begin{align}
\mathcal{B}_5(N)& =\int_{-\infty}^{+\infty}(L^{+}(x))^{4}\Theta(x)e(-Nx)\mathrm{d}x  \notag\\
&=\left(\intop_{|x|> K}+\intop_{\tau \le|x|\le K}+\intop_{|x|<\tau}\right)(L^{+}(x))^{4}\Theta(x)e(-Nx)\mathrm{d}x \notag\\
&=:\mathcal{B}_{5}^{(1)}(N)+\mathcal{B}_{5}^{(2)}(N)+\mathcal{B}_{5}^{(3)}(N).
\end{align}

Next, we give the estimation of the six terms in (4.4) and (4.5) by giving the following three Propositions.
\\
{\bf Proposition 4.1}. {\it For $1<c<\frac{1787}{1502},$ we have}
\begin{align*}
\mathcal{B}_{1}^{(1)}(N),\,\, \mathcal{B}_{5}^{(1)}(N)\ll 1.
\end{align*}
{\bf Proof}.
By the trivial estimate $L^\pm(x)\ll X^{1+\varepsilon}$ and Lemma 2.1, we get
\begin{align*}
\mathcal{B}_{1}^{(1)}(N),\,\, \mathcal{B}_{5}^{(1)}(N)& \ll X^{4+\varepsilon}\left(\frac{r}{\Delta}\right)^{r}\int_{K}^{+\infty}\frac{1}{x^{r+1}}\mathrm{d}x  \ll X^{4+\varepsilon}\left(\frac{r}{\Delta K}\right)^{r}\ll 1.
\end{align*}
\\
{\bf Proposition 4.2}. {\it For $1<c<\frac{1787}{1502},$ we have}
\begin{align*}
\mathcal{B}_{1}^{(2)}(N),\,\, \mathcal{B}_{5}^{(2)}(N) \ll \Delta X^{4-c-\varepsilon}.
\end{align*}
{\bf Proof}. For $\mathcal{B}_{1}^{(2)}(N)$, by H\"{o}lder's inequality, Lemma 2.1 and Lemma 2.11 (iii), we obtain
\begin{align*}
\mathcal{B}_{1}^{(2)}(N)& \ll\intop_{\tau\leq|x|\leq K}\big|L^{-}(x)\big|\big|L^{+}(x)\big|^3\big|\Theta(x)\big|\mathrm{d}x \ll\Delta\intop_{\tau\leq|x|\leq K}\left|L^{-}(x)\right|\left|L^{+}(x)\right|^{3}\mathrm{d}x \notag\\
&\ll\Delta\bigg(\intop_{\tau\leq|x|\leq K}\big|L^{-}(x)\big|^{4}\mathrm{d}x\bigg)^{\frac{1}{4}}\bigg(\intop_{\tau\leq|x|\leq K}\big|L^{+}(x)\big|^{4}\mathrm{d}x\bigg)^{\frac{3}{4}}\ll\Delta X^{4-c-\varepsilon}.
\end{align*}
Similarly, we have
\begin{align*}
\mathcal{B}_{5}^{(2)}(N) \ll \Delta\intop_{\tau\leq|x|\leq  K}\left|L^{+}(x)\right|^{4}\mathrm{d}x\ll \Delta X^{4-c-\varepsilon}.
\end{align*}
\\
{\bf Proposition 4.3}. {\it For $1<c<\frac{1787}{1502},$ we have }
\begin{align*}
\mathcal{B}_{1}^{(3)}(N)&=\mathcal{M}^{-}\big(\mathcal{M}^{+}\big)^{3}\mathcal{I}+O\big(\Delta{X^{4-c}}L^{11-A}\big),\\
\mathcal{B}_{5}^{(3)}(N)&=\left(\mathcal{M}^+\right)^4\mathcal{I}+O\big(\Delta{X^{4-c}}L^{11-A}\big),
\end{align*}
where $\mathcal{I}$ is deined by (4.6).
\\
{\bf Proof}. We shall give the evaluation of $\mathcal{B}_{1}^{(3)}(N)$, while $\mathcal{B}_{5}^{(3)}(N)$ can be considered in the same way. We define
\begin{align}
\mathcal{I}_0&=\int_{-\tau}^{\tau}I^4(x)\Theta(x)e(-Nx)\mathrm{d}x,\,\,\,\,\,\,\,\,\,\,\,\mathcal{I}=\int_{-\infty}^{+\infty}I^{4}(x)\Theta(x)e(-Nx)\mathrm{d}x.
\end{align}
It is easy to see that
\begin{align}
\mathcal{B}_{1}^{(3)}(N) = \mathcal{M}^{-}\big(\mathcal{M}^{+}\big)^{3}\mathcal{I} + \mathcal{M}^{-}\big(\mathcal{M}^{+}\big)^{3}\left(\mathcal{I}_0-\mathcal{I}\right) + \left(\mathcal{B}_{1}^{(3)}(N) - \mathcal{M}^{-}\big(\mathcal{M}^{+}\big)^{3}\mathcal{I}_0\right).
\end{align}
Since $\left|\frac{\mathrm{d}}{\mathrm{d}t}(t^cx)\right|\gg|x|X^{c-1}$ for $t\in(\mu X,X]$, then on applying Titchmarsh [17, Lemma 4.2], we find that $I(x) \ll |x|^{-1}X^{1-c}$. The estimate (2.4) together with Lemma 2.1 yields that
\begin{align}
\mathcal{M}^{-}\big(\mathcal{M}^{+}\big)^{3}\left(\mathcal{I}_0-\mathcal{I}\right)\ll&\:{\mathcal{L}}^4\int_\tau^{+\infty}\left|I(x)\right|^4\big|\Theta(x)\big|\mathrm{d}x
\ll\Delta {\mathcal{L}}^4 X^{4-4c}\int_\tau^{+\infty}\frac{1}{x^4}\mathrm{d}x \notag\\
\ll&\:\Delta {\mathcal{L}}^4 X^{4-4c}\tau^{-3}\ll\Delta {\mathcal{L}}^4X^{4-c-3\xi}.
\end{align}
Moreover, we use the identity
\begin{align*}
&L^{-}(x)\big(L^{+}(x)\big)^{3}-\mathcal{M}^{-}\big(\mathcal{M}^{+}\big)^{3}I^{4}(x) \notag\\
=&\left(L^{+}(x)\right)^{3}\big(L^{-}(x)-\mathcal{M}^{-}I(x)\big)+\mathcal{M}^{-}I(x)\big(L^{+}(x)\big)^{2}\big(L^{+}(x)-\mathcal{M}^{+}I(x)\big) \notag\\
&+\mathcal{M}^{-}\mathcal{M}^{+}I^{2}(x)L^{+}(x)\Big(L^{+}(x)-\mathcal{M}^{+}I(x)\Big)+\mathcal{M}^{-}\Big(\mathcal{M}^{+}\Big)^{2}I^{3}(x)\Big(L^{+}(x)-\mathcal{M}^{+}I(x)\Big)\notag
\end{align*}
and the estimate (2.4) to find that
\begin{align*}
\left|L^-(x)\big(L^+(x)\big)^3-\mathcal{M}^-\big(\mathcal{M}^+\big)^3I^4(x)\right|\ll X{\mathcal{L}}^{3-A}\Big(\big|L^{+}(x)\big|^{3}+\left|I(x)\right|^{3}\Big).
\end{align*}
From the above inequality,  together with Lemma 2.1 and  Lemma  2.8 (ii),(iv) gives that
\begin{align}
\mathcal{B}_{1}^{(3)}(N)-\mathcal{M}^-\big(\mathcal{M}^+\big)^3\mathcal{I}_0 &\ll\Delta X{\mathcal{L}}^{3-A}\left(\intop_{|x|\le \tau}\big|L^+(x)\big|^3\mathrm{d}x+\intop_{|x|\le \tau}\big|I(x)\big|^3\mathrm{d}x\right)\notag\\
&\ll \Delta X^{4-c}{\mathcal{L}}^{11-A}.
\end{align}
Combining (4.7)-(4.9), we have
\begin{align*}
\mathcal{B}_{1}^{(3)}(N)=\mathcal{M}^{-}\big(\mathcal{M}^{+}\big)^{3}\mathcal{I}+O\big(\Delta{X^{4-c}}{\mathcal{L}}^{11-A}\big).
\end{align*}
We proceed with $\mathcal{B}_{5}^{(3)}(N)$ in the same way and prove that
\begin{align*}
\mathcal{B}_{5}^{(3)}(N)=\big(\mathcal{M}^+\big)^4\mathcal{I}+O\big(\Delta{X^{4-c}}{\mathcal{L}}^{11-A}\big).
\end{align*}
On concluding (4.2),(4.4)-(4.5) and Propositions 1-3, we have
\begin{align}
\mathcal{B}(N) \geq \big(4\mathcal{M}^{-}-3\mathcal{M}^{+}\big)\big(\mathcal{M}^{+}\big)^{3}\mathcal{I}+O\big(\Delta{X^{4-c}}{\mathcal{L}}^{11-A}\big).
\end{align}
From Lemmas 2.2 (ii) and (2.9) (ii), we obtain
\begin{align}
&\big(4\mathcal{M}^{-}-3\mathcal{M}^{+}\big)\big(\mathcal{M}^{+}\big)^{3}\mathcal{I} \notag\\
\geq\,\,&\left(4f\bigg(\frac{\log D}{\log z_2}\bigg)-3F\bigg(\frac{\log D}{\log z_2}\bigg)\bigg)\big(1+O({\mathcal L}^{-1/3})\big)\mathfrak{P}^4\mathcal{I}\right.  \notag\\
=\,\,&\left(4f\bigg(\frac{62451}{20000}\bigg)-3F\bigg(\frac{62451}{20000}\bigg)\bigg)\mathfrak{P}^4\mathcal{I}+O\big(\Delta{X^{4-c}}{\mathcal{L}}^{-\frac{13}{3}}\big)\right.  \notag\\
=\,\,&\frac{160000e^\gamma}{62451}\bigg(\log\frac{42451}{20000}-\frac34\bigg(1+\int_2^{\frac{42451}{20000}}\frac{\log(t-1)}t\mathrm{d}t\bigg)\bigg)\mathfrak{P}^4\mathcal{I}
+O\big(\Delta{X^{4-c}}{\mathcal{L}}^{-\frac{13}{3}}\big) \notag\\
\ge\,\,&0.00027\times\mathfrak{P}^{4}\mathcal{I}+O\big(\Delta{X^{4-c}}{\mathcal{L}}^{-\frac{13}{3}}\big).
\end{align}
According to (4.10)-(4.11) and Lemma 2.9 (ii), we deduce that
\begin{align*}
\mathcal{B}(N) \gg\Delta X^{4-c}{\mathcal{L}}^{-4}.
\end{align*}
Therefore, $\mathcal{B}(N)>0$ and the inequality (1.4) has a solution in primes $p_1,p_2,p_3,p_4$ such that each of the numbers $p_{i}+2\,(i=1,2,3,4)$ had at most [$\eta_2^{-1}]=\big[\frac{93801402}{35740000-30040000c}\big]$  prime factors. This completes the proof of the Theorem 2.

{\bf{Acknowledgement}}. The author would like to thank the anonymous referee for his/her patience and time in refereeing this manuscript.

\end{document}